\pdfoutput=1 
\documentclass[12pt]{article} 
\usepackage[utf8]{inputenc} 
\usepackage[english]{babel} 
\usepackage{times} 

\usepackage{amsmath}
\usepackage{booktabs}
\usepackage{graphicx}
\usepackage{geometry}
\usepackage{gbt7714}
\usepackage{pdfpages}
\usepackage{type1cm}
\usepackage{footnote}
\usepackage{listings}
\usepackage{xcolor}
\usepackage{amssymb}
\usepackage{amsfonts}
\usepackage{multirow}
\usepackage{subfigure}   
\usepackage{amsthm}
\usepackage{setspace}
\usepackage{subcaption}
\usepackage{etoolbox}
\usepackage{hyperref}  
\usepackage{thumbpdf,lmodern}
\usepackage{framed}
\usepackage{fancyheadings}
\usepackage{rotating}

\newtheorem{definition}{Definition} 
\newtheorem{theorem}{Theorem}
\newtheorem{lemma}{Lemma}

\numberwithin{equation}{section} 

\newcommand{\upcite}[1]{\textsuperscript{\cite{#1}}}

\begin{document}
\pagestyle{fancy}

\title{Trimmed Mean for Partially Observed Functional Data}
\author{
	Yixiao Wang\thanks{Graduated from the University of Science and Technology of China (USTC).} \\
	\texttt{yixiao.wang@duke.edu}
}

\date{}
\maketitle

\begin{abstract}
	In recent years, partially observable functional data has gained significant attention in practical applications and has become the focus of increasing interest in the literature. In this thesis, we build upon the concept of data integration depth for partially observable functions, as proposed by Elias et al. (2023) \upcite{EliasEtAl2023}, and the trimmed mean estimator method along with its consistency proof introduced by Fraiman and Muniz (2001) \upcite{FraimanMuniz2001} for completely observable functions. We introduce the concept of trimmed mean specifically for partially observable functional data. Additionally, we address several theoretical and practical issues, including a proof of the strong consistency of the proposed trimmed mean, and we provide a simulation study. The results demonstrate that our estimator outperforms the ordinary mean in terms of accuracy and robustness when applied to partially observable functional data.
\end{abstract}

\textbf{Keywords:} Trimmed mean; Partially observed functional data; Functional depth; Robust

\fancyfoot{}
\cfoot{\thepage} 

\newpage

\section{Introduction}

\subsection{Background and Literature Review}

Over the past 80 years, functional data analysis (FDA) has emerged as one of the main branches of modern statistics. Analyzing functional data directly is often more effective and concise compared to traditional data-based approaches, particularly when dealing with dynamic data streams or large datasets. To better understand the process of dynamic data collection, Ramsay (1982) \upcite{Ramsay1982} highlighted that advanced modern data collection systems can now capture a series of functional observations. Ramsay emphasized that this data collection process should be conceptualized as dynamic rather than static, contrasting with the static datasets traditionally studied. Consequently, relying on conventional data analysis methods may result in information loss or biased model estimations.

From another perspective, when the volume of data is large, Ramsay and Silverman (1997) \upcite{RamsaySilverman1997} demonstrated that analyzing functional data is often more suitable than handling large finite-dimensional vectors derived from discrete approximations of functions. This is particularly relevant for computational reasons, as efficiently computing multivariate depth in high-dimensional spaces is nearly infeasible. This underscores the necessity of adapting traditional statistical methods to accommodate functional data. Hence, extending these methods is critical for preserving the richness and structure inherent in functional datasets.

Another important traditional statistical concept is order statistics, which concerns the properties of data rankings. This concept has broad applications, ranging from basic measures like medians and quartiles to various robust estimation techniques. In the univariate case, order statistics are straightforward, based on the ranking of real numbers, and their applications are well-documented in the literature. For example, L-estimates, defined as linear combinations of order statistics, are widely used as robust location estimators. These include specific cases such as the median (derived from a single point in an odd sample) and the mean (using all points). Trimmed means, a particular type of L-estimate, are calculated as the average of the central $(1-\alpha) n$ observations, where $(0 \leq \alpha < 1)$ varies. This approach balances efficiency and robustness by utilizing more data than the median while mitigating the influence of outliers. In higher dimensions, the concept of order statistics becomes more complex. Definitions proposed by Tukey (1975) \upcite{Tukey1975} and Fraiman and Meloche (1999) \upcite{FraimanMeloche1999} rely on various notions of data depth, which measures the "centrality" of multivariate data points within a given dataset to extend the concept of ranking in high dimensional case. Although these definitions are differ greatly for multivariate data, they are very similar when applied to univariate data.

Extending the concept of infinite-dimensional spaces, particularly for functional data, is both important and natural. Starting with the seminal work of Fraiman and Muniz (2001) \upcite{FraimanMuniz2001}, various definitions of depth for functional data have been proposed. The theoretical properties of these functional depths have been extensively studied (e.g., Nieto-Reyes and Battey (2016) \upcite{NietoReyesBattey2016}). A central idea in functional depth lies in integrating pointwise depths across the domain of the function. Functional depth has proven to be a powerful tool for numerous applications, such as visualizing functional data (Hyndman and Shang (2010) \upcite{HyndmanShang2010}), detecting functional outliers (Arribas-Gil and Romo (2014) \upcite{ArribasGilRomo2014}), and classifying functional data (Li, Cuesta-Albertos, and Liu (2012) \upcite{LiEtAl2012}).

However, in real-world scenarios, functional data are often incomplete due to issues like missing recordings, equipment failures, or other practical constraints. Extending the concept of depth to handle such incomplete functions is a challenging but fascinating problem. Partially observed functional data are prevalent across various fields of research. For instance, in medical studies, missing observations may arise from patients skipping medical visits or equipment failing to record data (James, Hastie, and Sugar (2000) \upcite{JamesEtAl2000}). Similarly, in electricity markets, supply functions are incomplete because suppliers and buyers negotiate prices and quantities based on changing market conditions (Kneip and Liebl (2020) \upcite{KneipLiebl2020}). As a result, the analysis of partially observed functional data has gained significant attention in recent years.Mathematically, partially observed functional data can be defined as follows: Consider random functions sampled on a compact set $[a, b]$ to $\mathbb{R}$. In cases where the recorded data for these functions are available only on a compact subset of $[a, b]$, each functional observation becomes partially observable. This formulation underpins the challenges and opportunities for analysis in such settings.

For partially observed functional data, Elias et al. (2023) \upcite{EliasEtAl2023} introduced a robust depth function called Partially Observed Integrated Functional Depth (POIFD). This function extends the toolkit for analyzing complex functional data settings, including sparse cases, and facilitates the development of robust statistics. POIFD provides critical support for other techniques in the context of partially observed functional data. Building on this foundational work, we propose a natural $\alpha$-trimmed mean concept for partially observed functional data. Specifically, the $\alpha$-trimmed mean is defined as the average of the most central (i.e., deepest) $1-\alpha$ proportion of observations based on the integral depth of partially observed functions. Compared to the ordinary mean, our approach demonstrates superior efficiency and robustness, particularly for incomplete or noisy data.

\subsection{Organization of the Paper}

In Chapter 2, we review the concept of functional data depth. In Chapter 3, we define the partially observable functional $\alpha-$trimmed mean based on the partially observable functional depth and propose the strong consistency property of trimmed mean. In Chapter 4, we present the results of a simulation study that compares the performance of several location estimators. All proofs, figures, and tables are included in the appendix, and the corresponding code is available at \href{https://github.com/Yixiao-Wang-Stats/TMoPOFD}{GitHub}.

\section{Functional Data Depth}
\subsection{Data Depth}
A data depth is a function to indicate how deep or "centrality" of a data point within a given data cloud or a given probability distribution. Although these definitions vary significantly for data in high dimensions, they become quite similar when applied to data in one dimension. We now briefly describe two of them.

Let $X_{1}, X_{2}, \ldots, X_{n}$ be independent and identically distributed (i.i.d.) random vectors in $\mathbb{R}^{d}$, and the common distribution $F_{X}$.

\textbf{Tukey's Depth.} After given the set of n points $\mathcal {X}_{n} = \{x_{1},\dots,x_{n}\}$, which is an observation of $\mathcal \{X_{1}, X_{2}, \ldots, X_{n}\}$ in d-dimensional space,  the sample version of Tukey's depth of a point $x$ is the smallest fraction (or number) of points in any closed halfspace that contains $x$.The Tukey's depth at $x$ is defined as
$$ 
T D_{n}(x;{\mathcal {X}}_{n})=\inf _{v\in \mathbb {R} ^{d},\|v\|=1}{\frac {1}{n}}\sum _{i=1}^{n}\mathbf {1} \{v^{T}(x_{i}-x)\geq 0\},
$$
and the population Tukey's depth of x with regret to the common distribution $F_{X}$ is
$$
T D(x;F_{X})=\inf _{v\in \mathbb {R} ^{d},\|v\|=1}P(v^{T}(X-x)\geq 0).
$$

In the univariate case, we assume $n$ is large enough for us to approximate $F$ with $F_n$, then we have

$$
T D(x;F_{X})=\min \left\{F(x), 1-F\left(x^{-}\right)\right\} \quad \text { and } \quad T D_{n}(x;{\mathcal {X}}_{n})=\min \left\{F_{n}(x), 1-F_{n}\left(x^{-}\right)\right\} \text {. }\nonumber
$$

\textbf{Simplicial Depth.} The simplicial depth of a point $x$ in d-dimensional Euclidean space, with respect to the set of n points $\mathcal {X}_{n} = \{x_{1},\dots,x_{n}\}$, which is an observation of $\{X_{1}, X_{2}, \ldots, X_{n}\}$ in d-dimensional space, is the number of d-dimensional simplices (the convex hulls of sets of $d+1$ sample points) that contain $x$. The Tukey's depth at $x$ is defined as

$$
S D_{n}(x; \mathcal {X}_{n})=\frac{1}{\tbinom {n}{d+1}}\#\{x\in S\left[x_{1}, x_{2}, \ldots, x_{n}\right]\} ,
$$

The population Simplicial depth at $x$ is defined as

$$
S D(x; F_{X})=P\left(x\in S\left[X_{1}, X_{2}, \ldots, X_{n}\right]\right) ,
$$

In the univariate case, we assume $n$ is large enough for us to approximate $F$ with $F_n$, then we have

$$
S D(x;F_X)=2 F(x)\left(1-F\left(x^{-}\right)\right),\quad \text { and } \quad S D_{n}(x; \mathcal {X}_{n} )=2 F_{n}(x)\left(1-F_{n}\left(x^{-}\right)\right) \text {. }\nonumber
$$

\subsection{Integrated Functional Depth}

Consider defining a depth of functional data in a finite-dimensional and univariate case. Let $\mathcal{P}$ denote the set of all probability measures on $(\mathbb{R}, \mathcal{B}(\mathbb{R}))$, where $\mathcal{B}(\mathbb{R})$ is the $\sigma$-algebra on $\mathbb{R}$.

Without loss of generality, we consider a functional data defined on the interval $[0,1]$. Specifically, in what follows, $X: [0,1] \rightarrow \mathbb{R}$ is a stochastic process with continuous paths, $P$ is the probability distribution of $X$, and $P_{t}$ is the marginal distribution of $X(t)$.

Considering a function $\mathrm{D}: \mathbb{R} \times \mathcal{P} \rightarrow [0,1]$ to be with some kind of depth. We assume that this function satisfies properties $\mathrm{D}{1}$ to $\mathrm{D}{7}$ as proposed by Nagy et al. (2016)\upcite{NagyEtAl2016}. Using this depth function $\mathrm{D}$ we can define a functional depth for functional data. Specifically, we recall the definition of the Integrated Functional Depth (Claeskens et al., 2014\upcite{ClaeskensEtAl2014}; Nagy et al., 2016\upcite{NagyEtAl2016}).

\begin{definition}(Integrated Functional Depth (IFD))
	Given a continuous function $x: [0,1] \rightarrow \mathbb{R}$, a univariate depth $\mathrm{D}$, and a weight function $w: [0,1] \rightarrow [0, \infty)$ satisfying $\int_{0}^{1} w(t) dt = 1$, the Integrated Functional Depth of $x$ with respect to $P$ is defined as $\int_{0}^{1} w(t) d t=1$, 
	
	$$
	\operatorname{IFD}_{w}(x, P)=\int_{0}^{1} \mathrm{D}\left(x(t), P_{t}\right) w(t) d t.
	$$
\end{definition}

\subsection{Partially Observed Integrated Functional Depth (POIFD)}
\subsubsection{Basic Setting}
Let $X_{1}, \ldots, X_{n}$ be $n$ independent realizations of $X$. We consider the scenario where the realizations $X_{1}, \ldots, X_{n}$ are only partially observed. To construct this setting, similarly to Delaigle and Hall (2013)\upcite{DelaigleHall2013}, we consider a random observation mechanism $Q$ that generates a compact set on $[0,1]$, which represents the domain over which the functional data are observed. Specifically, we assume that the compact set generated by $Q$ consists of a finite collection of closed intervals with strictly positive Lebesgue measure. Let $\mathcal{O}$ be a set generated by the mechanism $Q$, and let $\mathcal{O}_{1}, \ldots, \mathcal{O}_{n}$ be independent realizations of $\mathcal{O}$. Then, for $1 \leq i \leq n$, the functional data $X_{i}$ are only observed on $\mathcal{O}_{i}$.

We assume that the probability measure $P$ and the mechanism $Q$ are independent, and that $\left(X_{1}, \mathcal{O}_{1}\right), \ldots, \left(X{n}, \mathcal{O}_{n}\right)$ are independent and identically distributed realizations from $P \times Q$. This assumption is known as Missing-Completely-at-Random (MCAR), which is a standard assumption in the literature on partially observed functional data (e.g., Kraus 2015\upcite{Kraus2015}; Kneip and Liebl 2020\upcite{KneipLiebl2020}).

\subsubsection{Population Form}
Since $X$ is only observed on the set $\mathcal{O}$, we define its depth by restricting an integrated functional depth to the set $\mathcal{O}$. For $t \in [0,1]$, to achieve this, let $Q(t) = \mathbb{P}(\mathcal{O} \ni t)$ denote the probability that the random set $\mathcal{O}$ covers the point $t$. Without loss of generality, we assume that $Q(t) > 0$ for any $t \in [0,1]$. Furthermore, we consider a positive bounded continuous function $\phi$ defined on $[0,1]$; such a function $\phi$ could be, for example, the identity map on $[0,1]$. We then define the following weight function, restricted to the compact set $\mathcal{O}$:

$$
w_{\phi}(t \mid \mathcal{O})=\frac{\phi(Q(t))}{\int_{\mathcal{O}} \phi(Q(t)) d t}.
$$
We now recall the definition of the Partially Observed Integrated Functional Depth (POIFD) for any continuous function $x: [0,1] \rightarrow \mathbb{R}$ restricted to a compact set $\mathcal{O} \subset [0,1]$ (Elias et al. 2023\upcite{EliasEtAl2023}).

\begin{definition}(Partially Observed Integrated Functional Depth (POIFD))
	For $(x, \mathcal{O})$, the POIFD with respect to $P \times Q$ is defined as
	
	$$
	\operatorname{POIFD}((x, \mathcal{O}), P \times Q)=\int_{\mathcal{O}} \mathrm{D}\left(x(t), P_{t}\right) w_{\phi}(t \mid \mathcal{O}) d t,
	$$
	
	where $w_{\phi}(t \mid \mathcal{O})$ is the weight function defined above.
\end{definition}
It is easy to see that if the data are fully observed, i.e., if $\mathbb{P}(\mathcal{O} \ni t) = 1$ for any $t \in [0,1]$, then the proposed POIFD corresponds exactly to the IFD with the weight function $w(t) \equiv 1$.

\subsubsection{Sample Form}
We now recall the sample version of POIFD for computation on finite samples. Let $P_{n}$ denote the empirical distribution obtained by assigning a weight of $1/n$ to each sample curve $X_{1}, \ldots, X_{n}$. Similarly, let $Q_{n}$ denote the empirical distribution obtained by assigning a weight of $1/n$ to each set $\mathcal{O}_{1}, \ldots, \mathcal{O}_{n}$. Define $I(t) := \{1 \leq i \leq n: t \in \mathcal{O}_{i}\}$ and $q_{n}(t) = \frac{\# I(t)}{n}$. Finally, let $P_{t, n}$ represent the empirical distribution function of the univariate sample $\{X_{i}(t), i \in I(t)\}$.

We define the sample version of POIFD using a plug-in approach(Elias et al. 2023\upcite{EliasEtAl2023}):

\begin{definition}(Sample Version of POIFD)
	$$
	\begin{aligned}
		\operatorname{POIFD_{n}}\left((x, \mathcal{O}), P_{n} \times Q_{n}\right) =\frac{\int_{\mathcal{O}} \mathrm{D}\left(x(t), P_{t, n}\right) \phi\left(q_{n}(t)\right) d t }{\int_{\mathcal{O}} \phi\left(q_{n}(t)\right) d t}.
	\end{aligned}
	$$
\end{definition}

In practical applications, we only need to obtain a discrete version of the functional data $X_{i}$ on a discrete evaluation grid $\{t_{1}, \ldots, t_{T}\}$, where $0=t_{1}<t_{2}<\cdots<t_{T}=1$ (often selected equidistantly for simplicity). Note that for many $t_{\ell}$, $\ell=1, \ldots, T$, such evaluations may indeed be missing, which is consistent with the partially observed nature of the data considered in this work. We can then define the sample version of the POIFD using the standard Riemann approximation.

\begin{definition}
	(Sample Partially Observed Integrated Functional Depth (Sample POIFD)) Given partially observed functional data $\left(X_{1}, \mathcal{O}_{1}\right), \ldots, \left(X_{n}, \mathcal{O}_{n}\right)$, evaluated on a common grid $\{t_{1}, \ldots, t_{T}\}$, where $0=t_{1}<t_{2}<\cdots<t_{T}=1$, the sample version of the Partially Observed Integrated Functional Depth is defined as
	You can complete the definition by specifying the formula that represents the sample POIFD, using the given discrete evaluation grid.
	
	$$
	\begin{aligned}
		\operatorname{POIFD}_{T}\left((x, \mathcal{O}), P_{n} \times Q_{n}\right) =\sum_{t_{\ell} \in \mathcal{O}} \mathrm{D}\left(x\left(t_{\ell}\right), P_{t_{\ell}, n}\right) \frac{\phi\left(q_{n}\left(t_{\ell}\right)\right)}{\sum_{t_{\ell} \in \mathcal{O}} \phi\left(q_{n}\left(t_{\ell}\right)\right)}.
	\end{aligned}
	$$
\end{definition}

\section{The Trimmed Mean for Partially Observed Functional Data and Related Theorems}
\subsection{Definition}

The trimmed mean for partially observed functions is defined as the mean of the $n - [n\alpha]$ deepest observed values, where $\alpha$ is the ratio of trimming.

More precisely, for $\beta > 0$, we define the sample version of trimmed mean for partially observed functional data as

$$
\hat{\mu}_{n}(t)=\frac{\sum_{i=1}^{n} \mathbf{1}_{[\beta,+\infty)}\left(POIFD_{n}\left(X_{i}\right)\right)\mathbf{1}_{X_{i}(t)\text{ is observed}} X_{i}(t)}{\sum_{i=1}^{n} \mathbf{1}_{[\beta,+\infty)}\left(POIFD_{n}\left(X_{i}\right)\right)\mathbf{1}_{X_{i}(t)\text{ is observed}}} ,\quad t\in\left[0,1\right],
$$

where $\beta$ satisfies

$$
\frac{1}{n} \sum_{i=1}^{n} \mathbf{1}_{[\beta,+\infty)}\left(POIFD_{n}\left(X_{i}\right)\right) \approx 1-\alpha .
$$

Here, $\mathbf{1}_{A}$ denotes the indicator function of the set $A$.

Similarly, we define the population version of the trimmed mean for partially observed functions as

$$
\mu(t)=\frac{E\left[\mathbf{1}_{\left[\beta,+\infty\right)}\left(POIFD\left(X\right)\right)\mathbf{1}_{X(t)\text{ is observed}} X(t)\right]}{{E\left[\mathbf{1}_{\left[\beta,+\infty\right)}\left(POIFD\left(X\right)\right)\mathbf{1}_{X(t)\text{ is observed}}\right]}} , \quad t\in\left[0,1\right],
$$

We also define 

$$
\hat{\mu}(t)=\frac{\sum_{i=1}^{n} \mathbf{1}_{[\beta,+\infty)}\left(POIFD\left(X_{i}\right)\right)\mathbf{1}_{X_{i}(t)\text{ is observed}} X_{i}(t)}{\sum_{i=1}^{n} \mathbf{1}_{[\beta,+\infty)}\left(POIFD\left(X_{i}\right)\right)\mathbf{1}_{X_{i}(t)\text{ is observed}}} , \quad t\in\left[0,1\right],
$$

for our proof process.

\subsection{Strong Consistency Results}

Next, we will demonstrate the strong consistency of the empirical functional depth $POIFD_{n}$ with respect to its population counterpart $POIFD$ over an appropriate function space, and derive the strong consistency of the trimmed mean estimator. This part is an extension of trimmed mean from functional data to partially observed functional data based on the Fraiman and Muniz (2001) \upcite{FraimanMuniz2001}'s work. We recall the following two assumptions setting:

H1. For a sufficiently large constant $A$, let $$\operatorname{Lip}[0,1]=\{x:[0,1] \rightarrow R, \text{x is a Lipschitz function with a Lipschitz constant less than or equal to A}\},$$
which is the function space in which the random process $X_{1}(t)$ takes its values.

H2. There exists a constant $c > 0$ such that
$$
E\left(\lambda\left(\left\{t: X_{1}(t) \in[u(t), u(t)+c \epsilon]\right\}\right)\right)<\epsilon / 2.
$$
Here, $\lambda$ represents the Lebesgue measure on $\mathbb{R}$, and $u \in \operatorname{Lip}[0,1]$.

We now present the strong consistency results for the trimmed mean of partially observed functions.
\subsubsection{The Case of a Fixed Weight Function}
We first consider a given weight function $\omega(t)$, where $\sup_{t \in [0,1]} |\omega(t)| \leq C$ for some constant $C$. Under these conditions, we can prove the following theorem:

\begin{theorem}
	Under the conditions H1 and H2, if the weight function is $\omega(t)$ and $\sup_{t \in [0,1]} |\omega(t)| \leq C$ for some constant $C$, let
	$$
	K_{n}(x)=\int_{0}^{1} F_{n, t}(x(t))\omega(t) d t, \quad \text { and } \quad K(x)=\int_{0}^{1} F_{t}(x(t))\omega(t) d t.
	$$
	where $F_{n,t}(x(t))$ and $F_{t}(x(t))$ are the empirical and population distribution functions, respectively. We have
	$$
	\sup _{x \in \operatorname{Lip}[0,1]}\left|K_{n}(x)-K(x)\right|\xrightarrow{\text{a.s.}}0.
	$$
	and
	$$
	\sup _{x \in \operatorname{Lip}[0,1]}\left|POIFD_{n}(x)-POIFD(x)\right|\xrightarrow{\text{a.s.}}0 .  
	$$
\end{theorem}
\subsubsection{The Case of Weighted Functions for Partially Observed Functions}
Next, we consider the weighted functions used in this paper, with the corresponding population and empirical versions defined as:
$$
w_{\phi}(t \mid \mathcal{O})=\frac{\phi(Q(t))}{\int_{\mathcal{O}} \phi(Q(t)) d t}\quad \text{ and }\quad w_{\phi,n}(t \mid \mathcal{O})=\frac{\phi(q_n(t))}{\int_{\mathcal{O}} \phi(q_n(t)) d t}.
$$
\begin{theorem}
	Under the conditions H1 and H2, let
	$$
	K_{n}(x)=\int_{\mathcal{O}} F_{n, t}(x(t)) w_{\phi, n}(t\mid \mathcal{O}) d t, \quad \text { and } \quad K(x)=\int_{\mathcal{O}} F_{t}(x(t)) w_{\phi}(t \mid \mathcal{O}) d t.
	$$
	where $F_{n,t}(x(t))$ and $F_{t}(x(t))$ are the empirical and population distribution functions, respectively. Then we have
	$$
	 \sup_{x \in \operatorname{Lip}[0,1]} \left| K_{n}(x) - K(x) \right| \xrightarrow{\text{a.s.}} 0,
	$$
	and
	$$
	 \sup_{x \in \operatorname{Lip}[0,1]} \left| POIFD_{n}(x) - POIFD(x) \right| \xrightarrow{\text{a.s.}} 0,
	$$
	where $POIFD_{n}(x)$ and $POIFD(x)$ are the sample and population partially observed integrated functional depths, respectively.
\end{theorem}
\begin{theorem}
	If the stochastic process $X_{1}(t)$ takes values in an arbitrary
	space $E[0,1] := E$, where
	
	$$
	 \sup _{x \in E}\left|K_{n}(x)-K(x)\right|\xrightarrow{\text{a.s.}}0,  
	$$
	Then
	$$
	\hat{\mu}_{n} \rightarrow \mu \text { a.s., }
	$$
	where $\hat{\mu}_{n}$ is the empirical and population trimmed means, respectively. 
\end{theorem}

In particular, under the assumptions H1 and H2, we have $\hat{\mu}_{n}(t) \rightarrow \mu(t)$ almost surely, which proves the strong consistency of the trimmed mean for partially observed functional data. All the proofs can be found in the appendix, which are technical.

\section{Simulation Experiment}
\subsection{Model Setup}
In this section, we follow the simulation data assumptions from Fraiman and Muniz (2001)\upcite{FraimanMuniz2001}, and compare the trimmed mean estimates with the conventional mean under three contamination models.

The basic model (Model M1) consists of p functions, satisfying the following condition:

$$
X_{i}(t)=g(t)+e_{i}(t) \quad 1 \leq i \leq p,
$$
where $e_{i}(t)$ is a Gaussian stochastic process with zero mean, and the covariance function is given by
$$
E\left(e_{i}(t) e_{i}(s)\right)=\left(\frac{1}{2}\right)^{|t-s| p},
$$
Function $g(t) = 4t$ corresponds to the unpolluted model.

Next, we consider two types of contamination in the basic model: complete contamination and partial contamination (on the trajectory), while also accounting for both symmetric and asymmetric contamination.

In the case of symmetric complete contamination, model $M_2$ is given by the following formula:

$$
Y_{i}^{full}(t)=X_{i}(t)+\epsilon_{i} \sigma_{i} M \quad 1 \leq i \leq p,
$$

where $Y_{i}(t)$ is the partially observable functional data generated from $Y_{i}^{\text{full}}(t)$. Here, $\epsilon_{i}$ and $\sigma_{i}$ are independent sequences of random variables. $\epsilon_{i}$ takes the value 1 with probability $q$, and 0 with probability $(1-q)$, determining the proportion of contamination. It takes the values 1 and -1 with probability $1/2$ each, determining the direction of contamination. $M$ is the magnitude of the contamination (a constant).

In the case of asymmetric complete contamination, model $M_3$ is defined by the following formula:

$$
Y_{i}(t)=X_{i}(t)+\epsilon_{i} M \quad 1 \leq i \leq p,
$$

where $\epsilon_{i}$ and $M$ are as defined in model $M_2$.

For partial contamination, we consider model $M_4$ defined as follows:
$$
Y_{i}(t)=X_{i}(t)+\epsilon_{i} \sigma_{i} M \quad \text { for } t \geq T_{i} \quad 1 \leq i \leq p
$$

and

$$
Y_{i}(t)=X_{i}(t) \quad \text { for } t<T_{i},
$$
where $T_{i}$ is randomly chosen according to a uniform distribution on $(0,1)$.

We consider the following scenarios: for $p=50$ and $p=80$ curves, with an estimated number of selected grid points $len=200$, $q=0.1$, $M=5$ and $M=25$, and $\alpha=0.2$ and $0.3$, with observation proportions of 0.5 and 0.9. Each scenario was repeated $N=10$ times.

For constructing the partially observable functional data over random intervals: each function is observed over $m$ disjoint intervals, distributed along $[0,1]$, with a total expected observation proportion $p$. For this, for each function, a random sample of size $\lfloor(m-p)/p\rfloor$ is generated from a uniform distribution on $[0,1]$. We then consider the intervals $[u_{(i-1)}, u_{(i)}]$, where $u_{(i)}$ is the $i$-th order statistic of the sample. Thus, the function is observed over $m$ randomly chosen disjoint intervals. Note that the selection is random, but the $m$ intervals are guaranteed to be disjoint.

For simplicity in the simulation experiments, a center of 0.5 can be chosen. Specifically, each functional data point is observable over a domain from a randomly generated starting point to a randomly generated endpoint. The starting point is sampled from a uniform distribution on $[1/2-p, 1/2)$ if $p \leq 1/2$, or from a uniform distribution on $[0,1-p)$ if $p>1/2$. For the endpoint, if $p \leq 1/2$, it is sampled from a uniform distribution on $(1/2,1/2+p]$, or if $p>1/2$, from a uniform distribution on $(p,1]$. This results in a total (expected) observation proportion $p$.

\begin{figure}[htbp]
	\centering
	\subfigure[Partially observable functions]{
		\includegraphics[scale=0.7]{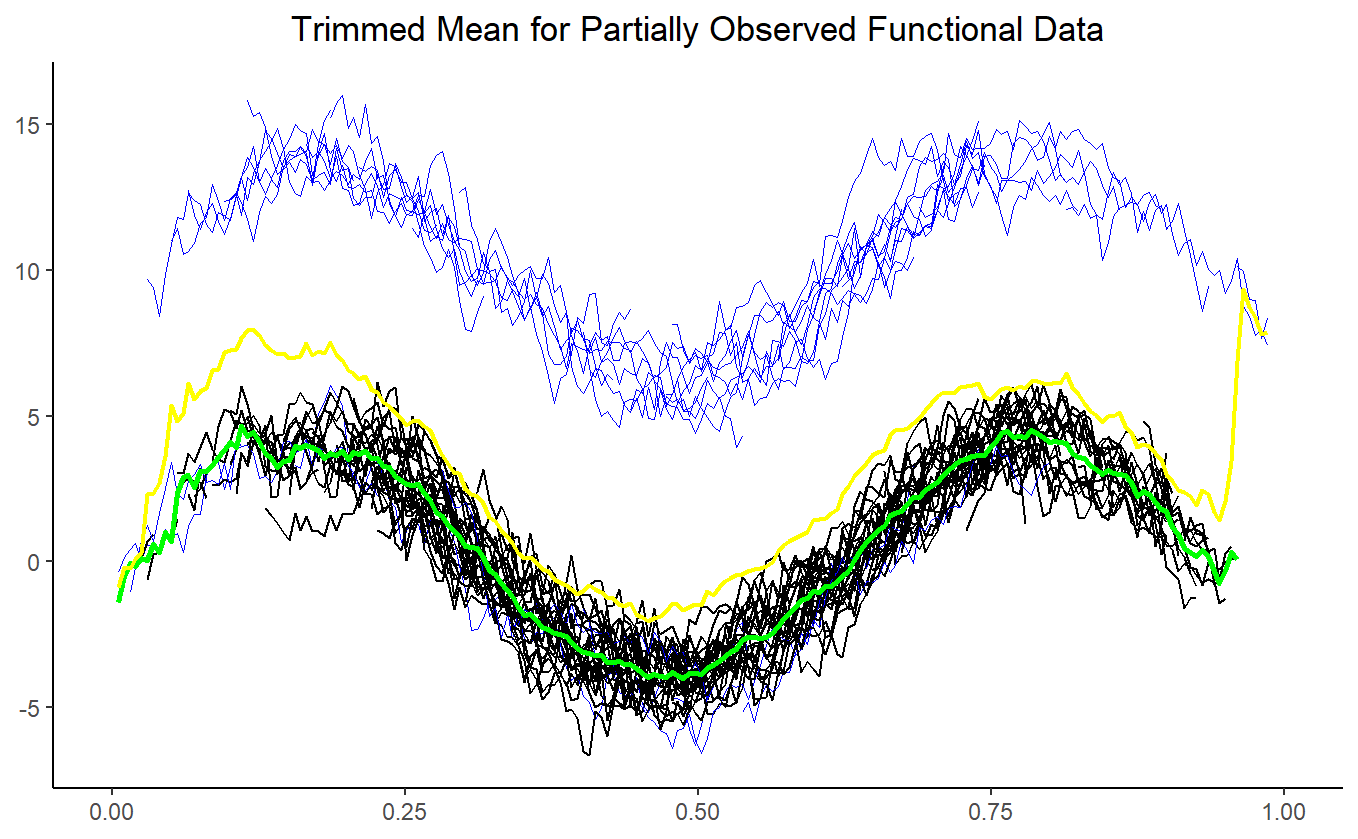}
	}
	\caption{Illustration of partially observable functions.}
	\label{figure0}
\end{figure}

As an example, see the Figure \ref{figure0}. The green line represents the trimmed mean, the yellow line represents the mean before trimming, and the blue lines represent the functions that were trimmed. We observe that all contamination can be identified and removed. In this figure, we use a sine function as the base setting, which slightly differs from our model setting.

\subsection{Experimental Data and Analysis}

In each case, for $p=50$ and $p=80$ curves, with an estimated number of selected grid points $len=200$, $q=0.1$, $M=5$ and $M=25$, and $\alpha=0.2$ and $0.3$, with observation proportions of 0.5 and 0.9, each scenario was repeated $N=10$ times.

Figures \ref{figure1}, \ref{figure2}, and \ref{figure3} (in the appendix) show the partially observable data generated under the same $P$, where $P = 50$, for the three types of contamination. Each figure has an upper panel showing the partially observed functions and a lower panel showing the proportion of observable functions $q_{n}$ at each point. The left side shows the complete data, while the right side shows the data after trimming with an observation proportion $\alpha = 0.3$. It can be seen that contamination data clearly lies away from the center of the original data on the left side of each figure. Meanwhile, on the right side of each figure, it is evident that the contamination values (i.e., the functions with the weakest centralization) are effectively trimmed off.

For each model, we consider mean and trimmed mean estimates.
$$
\hat{\mu}_{n}(t)=\frac{\sum_{i=1}^{n} \mathbf{1}_{X_{i}(t)\text{is observed}} X_{i}}{\sum_{i=1}^{n} \mathbf{1}_{X_{i}(t)\text{is observed}}} ,
$$

$$
\hat{\mu}_{n,\alpha}(t)=\frac{\sum_{i=1}^{n} \mathbf{1}_{[\beta,+\infty)}\left(POIFD_{n}\left(X_{i}\right)\right)\mathbf{1}_{X_{i}(t)\text{is observed}} X_{i}}{\sum_{i=1}^{n} \mathbf{1}_{[\beta,+\infty)}\left(POIFD_{n}\left(X_{i}\right)\right)\mathbf{1}_{X_{i}(t)\text{is observed}}} ,
$$

The trimming proportions $\alpha = 0.2$ and $0.3$ were selected.

For each of the 10 experiments, we evaluated the estimates at $\mathrm{I} = len$ equidistant points on $[0,1]$, and computed the integral error for each experiment:

$$
E I(j)=\frac{1}{I} \sum_{k=1}^{I}\left[f\left(\frac{k}{I}\right)-g\left(\frac{k}{I}\right)\right]^{2},
$$
where $f$ represents either $\hat{\mu}_{n}$ or $\hat{\mu}_{n, \alpha}$.

In the table, we report the average integral error of each estimate:

$$
E=\frac{1}{N} \sum_{j=1}^{N} E I(j)
$$

and its standard deviation:
$$
\left.s=\left(\frac{1}{N} \sum_{j=1}^{N}(E I(j)-E)^{2}\right)\right)^{1 / 2}.
$$
We also report a robustness measure for evaluating estimator performance:

$$
M=\operatorname{median}(E I(j)\text{, }j=1, \ldots N).
$$

We then averaged the error across the 10 experiments to produce the results in Table \ref{tab:tab1}, Table \ref{tab:tab2}, Table \ref{tab:tab3}, and Table \ref{tab:tab4}, which are shown in the appendix. All statistical results with the suffix "trim" refer to the trimmed values. It is evident that, in almost every setting, the trimmed mean proposed for partially observable functions significantly outperforms the ordinary mean in terms of average integral error, standard deviation, and robustness.

\newpage
\bibliography{refer} 
\newpage
\appendix
\section{Proof of Strong Consistency Results}

\subsection{Proof of Theorem 1}

\textbf{Proof:} We consider a slight modification of the proof of Theorem 3.1 in Fraiman and Muniz (2001) \upcite{FraimanMuniz2001}. Firstly we consider $Y(t)$ is also a realization of $X$. Thus $Y(t)$ still has the same distribution $F$ and marginal distribution $F_{t}$. Then by Fubini theorem, we have   
\begin{align}
	K(x)&=\int_{0}^{1} F_{t}(x(t))\omega(t) d t\nonumber\\
	&=\int_{0}^{1} P(Y(t)\leq x(t))\omega(t) d t\nonumber\\
	&=\int_{0}^{1} \int_{\Omega} \mathbf{1}_{\left(-\infty, x(t)\right]}(Y(t))\omega(t) dP dt\nonumber\\
	&=  \int_{\Omega}\int_{0}^{1} \mathbf{1}_{\left(-\infty, x(t)\right]}(Y(t)) \omega(t)dt dP \nonumber\\
	&:=\int_{\Omega}g_{x} dP\nonumber\\
	&:=Pg_{x}\nonumber,
\end{align}
where $g_x : L^1([0,1], \mathbb{R}) \to \mathbb{R}, \quad z(t) \mapsto \int_0^1 \mathbf{1}_{(-\infty, x(t)]}(z(t)) \omega(t) dt$, and $Pf$ is defined as $Pf = \int f dP$.

We define $\mathcal{F} := \{g_x; x\in Lip[0,1]\}$, which with envelope $F$ $\equiv$ $C$ trivially. By the Theorem  in Pollard (1984)\upcite{pollard1984}(Theorem 24 pag. 25), if we have $\log N_1(\epsilon, P_n, \mathcal{F}) = o_P(n)$, where $\log N_1(\epsilon, P_n, \mathcal{F})$ is the family’s entropy of $\mathcal{F}$ (Definition also can be found in Pollard (1984)\upcite{pollard1984} 23 pag. 25). Then we can prove
$$
\sup_{x \in \operatorname{Lip}[0,1]} \left| P_{n}g_x - Pg_{x} \right| \xrightarrow{\text{a.s.}} 0 ,
$$ 
which is 
$$
\sup_{x \in \operatorname{Lip}[0,1]} \left| K_{n}(x) - K(x) \right| \xrightarrow{\text{a.s.}} 0.
$$ 

In order to estimate $N_1(\epsilon, P_n, \mathcal{F})$, we considering for $x$, how to find $x^{\prime}$ to construct finite functions $g_{x^{\prime}}$, such that for all $g_{x}$, we can always find a $g_{x^{\prime}}$, which make $\left\|g_{x} - g_{x^{\prime}}\right\|_{L^{1}\left(P_{n}\right)}$ to be small enough.

More precisely, we define $y(t) = \min\{x(t),x^{\prime}(t)\}$ and $z(t) = \max\{x(t),x^{\prime}(t)\}$, then notice that
$$
\begin{aligned}
	\left\|g_{x} - g_{x^{\prime}}\right\|_{L^{1}\left(P_{n}\right)} &= \frac{1}{n} \sum_{i=1}^{n} \left|g_{x}\left(X_{i}\right) - g_{x^{\prime}}\left(X_{i}\right)\right| \\
	&= \frac{1}{n} \sum_{i=1}^{n} \left|\int_{0}^{1} \mathbf{1}_{(-\infty, x(t)]}\left(X_{i}(t)\right) w(t) \, dt - \int_{0}^{1} \mathbf{1}_{\left(-\infty, x^{\prime}(t)\right]}\left(X_{i}(t)\right) w(t) \, dt\right| \\
	&= \frac{1}{n} \sum_{i=1}^{n} \int_{0}^{1} \mathbf{1}_{[y(t), z(t)]}\left(X_{i}(t)\right) w(t) \, dt \\
	&\leq C \frac{1}{n} \sum_{i=1}^{n} \lambda\left\{t : X_{i}(t) \in [y(t), z(t)]\right\}.
\end{aligned}
$$

Note that if $\omega$ is dependent on $n$, as long as all the $\omega_{n}(t)$ have a common bound $C$, the derivation remains valid, which will be used in the proof of Theorem 2.

The remaining derivation is entirely consistent with Theorem 3.1 (pag. 11) of Fraiman and Muniz (2001) \upcite{FraimanMuniz2001}, in their Theorem 3.1 it has proven the situation without the constant C, the difference, $\frac{1}{n} \sum_{i=1}^{n} \lambda\left\{t : X_{i}(t) \in [y(t), z(t)]\right\}$, can be as small as we want, so that we can construct the the finite functions, then bounded the $N_1(\epsilon, P_n, \mathcal{F})$. A constant is not effect anything. 

\subsection{Proof of Theorem 2}

\begin{lemma}
	If \( X_{n} \xrightarrow{\text{a.s.}} X \) and \( Y_{n} \xrightarrow{\text{a.s.}} Y \), then \( \frac{X_{n}}{Y_{n}} \xrightarrow{\text{a.s.}}\frac{X}{Y} \). 
\end{lemma}

\textbf{Proof:} Let \( f : S \subset \mathbb{R}^2 \to \mathbb{R} \) be a continuous function, where \( S \) is the domain of \( f \). Assume that the random variables \( \{(X_n, Y_n), (X, Y)\} \) take values in \( S \) and that \( X_n \xrightarrow{\text{a.s.}} X \) and \( Y_n \xrightarrow{\text{a.s.}} Y \) Then \( f(X_n, Y_n) \xrightarrow{\text{a.s.}} f(X, Y) \).

To prove this assertion, let \( N = \{X_n \nrightarrow X\} \) and \( M = \{Y_n \nrightarrow Y\} \). Note that \( P(N \cup M) = 0 \). By the assumption of continuity, on \( (N \cup M)^c \), we have \( f(X_n, Y_n) \to f(X, Y) \), and \( P((N \cup M)^c) = 1 \). Therefore, when \( Y \) takes non-zero values, Lemma 1 holds.

\begin{lemma}
	(Egoroff's Theorem) Let \( (X, \mathfrak{a}, m) \) be a finite measure space, i.e., \( m(E) < \infty \), and let \( \{f_n\}, f \) be measurable functions that are almost everywhere real-valued. If \( f_n(x) \to f(x) \) a.e. for all \( x \in X \), then for any given \( \delta > 0 \), there exists a measurable subset \( E_\delta \) such that \( m(E_\delta) < \delta \), and \( f_n \to f \) uniformly on \( X \setminus E_\delta \).
\end{lemma}

\begin{lemma}
	If the real-valued functions \( f_{n}(t) \rightarrow f(t) \) a.e. on \( [0,1] \) and both \( f_n \) and \( f \) are bounded, then
	\[
	\int_{0}^{1}|f_{n}(t)-f(t)| dt \rightarrow 0.
	\]
\end{lemma}

\textbf{Proof:} By Egoroff's Theorem. For any given \( \delta, \epsilon > 0 \), there exists a measurable subset \( E_\delta \) such that \( m(E_\delta) < \delta \), and \( f_n \to f \) uniformly on \( [0,1] \setminus E_\delta \). This implies that there exists a sufficiently large \( n \) such that \( \sup_{t \in [0,1] \setminus E_\delta}|f_{n}(t)-f(t)| \leq \epsilon \).

\begin{align}
	\int_{0}^{1}|f_{n}(t)-f(t)| dt &= \int_{E_\delta}|f_{n}(t)-f(t)| dt + \int_{[0,1] \setminus E_\delta}|f_{n}(t)-f(t)| dt \nonumber \\
	&\leq m(E_\delta) C + \sup_{t \in [0,1]}|f_{n}(t)-f(t)| \nonumber \\
	&\leq \delta C+ \epsilon \rightarrow 0, \nonumber
\end{align}
where C is the bound of both $f_{n}$ and $f$. This completes the proof.

\begin{lemma}
	Let \( f, f_{n} : \Omega \times [0,1] \rightarrow \mathbb{R} \). If both \( f_n \) and \( f \) are bounded, and if for all \( t \in [0,1] \), we have \( f_{n}(t) \rightarrow f(t) \) a.s., then as \( n \rightarrow \infty \), we have
	\[
	\int_{0}^{1}|f_{n}(t)-f(t)| dt  \xrightarrow{\text{a.s.}} 0.
	\]
\end{lemma}

\textbf{Proof:} By Lemma 3, we know that Lemma 4 is equivalent to $\exists \Omega_{0}$ such that $P(\Omega_{0}) = 1$ and for all $\omega \in \Omega_{0}$, $f_{n, \omega}(t)  \xrightarrow{\text{a.e.}} f_{\omega}(t)$ as $n$ tends to infinity. Note that we use$f_{n, \omega}(t) = f_{n}(\omega, t)$ for simplicity of writing.

We will prove by contradiction. Consider the converse statement: $\forall \tilde{\Omega}$ with $P(\tilde{\Omega}) = 1$, there exists $\omega_{0} \in \tilde{\Omega}$ such that $f_{n,\omega_{0}} \nrightarrow f_{\omega_{0}}$ as $n$ tends to infinity, almost everywhere.

This is equivalent to stating that $\exists I_{\omega_{0}}$ with $m(I_{\omega_{0}}) > 0$ such that for all $t \in I_{\omega_{0}}$, $f_{n,\omega_{0}} \nrightarrow f_{\omega_{0}}$.

Let
\[
A \subset \Omega
\]
be defined as
\[
A = \{\omega_{0} \mid \exists I_{\omega_{0}}, m(I_{\omega_{0}}) > 0, \forall t \in I_{\omega_{0}}, f_{n,\omega_{0}} \nrightarrow f_{\omega_{0}}\}.
\]

Then for any $\tilde{\Omega}$ with $P(\tilde{\Omega}) = 1$, we have $\tilde{\Omega} \cap A \neq \emptyset$, which means $P(A) > 0$. Otherwise taking $\tilde{\Omega} = A^c$, as $P(\tilde{\Omega}) = P(A^c)=1$ and $\tilde{\Omega} \cap A = A^c\cap A = \emptyset$, leading to a contradiction.

Define
\[
A_{n} = \{\omega_{0} \mid \exists I_{\omega_{0}}, m(I_{\omega_{0}}) \geq \frac{1}{n}, \forall t \in I_{\omega_{0}}, f_{n,\omega_{0}} \nrightarrow f_{\omega_{0}}\}.
\]

Thus, we have $A = \bigcup_{n=1}^{\infty} A_{n}$. Since $P(A) > 0$, there exists $n_{0}$ such that $P(A_{n_{0}}) > 0$.

Let
\[
\tilde{A} = \{(\omega, t) \mid \omega \in A, t \in I_{\omega}\}.
\]

In particular, let $\tilde{A}_{n_{0}} = \{(\omega, t) \mid \omega \in A_{n_{0}}, t \in I_{\omega}\}$.

We now show that $P \times m(\tilde{A}) > 0$. In fact,
\begin{align}
	P \times m(\tilde{A}_{n_{0}}) &= \int \mathbf{1}_{\tilde{A}_{n_{0}}}(\omega, t) \, d m \, d P \nonumber \\
	&= \int_{A_{n_{0}}} d P \int_{0}^{1} \mathbf{1}_{\tilde{A}_{n_{0}}}(\omega, t) \, d m \nonumber \\
	&\geq \frac{1}{n_{0}} \int_{A_{n_{0}}} d P \nonumber \\
	&> 0. \nonumber
\end{align}

Since $\tilde{A}_{n_{0}} \subset \tilde{A}$, it follows that $P \times m(\tilde{A}) > 0$. This means that there is a positive measure set where $f_{n}(\omega, t) \nrightarrow f(\omega, t)$.

However, we have $\forall t$, $f_{n}(t) \rightarrow f(t)$ almost surely. This implies $\forall t$, $\exists \Omega_{t}$ with $P(\Omega_{t}) = 1$ such that $f_{n}(\omega, t) \rightarrow f(\omega, t)$. Thus,
\[
P \times m(\{(\omega, t) \mid \omega \in \Omega_{t}, t\in\left[0,1\right]\}) = \int_{0}^{1} d m \int_{\Omega_{t}} d P = 1,
\]
which implies $P \times m(\{(\omega, t) \mid f_{n}(\omega, t) \nrightarrow f(\omega, t)\}) = 0$. This is a contradiction!

\begin{lemma}
	If $X_{n} \rightarrow X$ a.s. and $f$ is a continuous function, then $f(X_{n}) \rightarrow f(X)$ a.s.
\end{lemma}

\textbf{Proof:} By the continuity of $f$, for all $\omega$, if $X_{n}(\omega) \rightarrow X(\omega)$, then $f(X_{n}(\omega)) \rightarrow f(X(\omega))$. Since $X_{n} \rightarrow X$ a.s., i.e., convergence occurs on a set of probability 1, it follows that $f(X_{n}) \rightarrow f(X)$ a.s.

We now proceed with the formal proof of \textbf{Theorem 2}.

\textbf{Proof:} By definition,
\[
K_{n}(x) = \int_{\mathcal{O}} F_{n, t}(x(t)) w_{\phi,n}(t \mid \mathcal{O}) \, d t \quad \text{and} \quad K(x) = \int_{\mathcal{O}} F_{t}(x(t)) w_{\phi}(t \mid \mathcal{O}) \, d t.
\]

Then we consider
\begin{align}
	&\left| K_{n}(x) - K(x) \right| \label{eq:Kn_diff}\nonumber\\
	&= \left| \int_{\mathcal{O}} \left[ F_{n, t}(x(t)) w_{\phi,n}(t \mid \mathcal{O}) - F_{t}(x(t)) w_{\phi}(t \mid \mathcal{O}) \right] \, d t \right| \nonumber \\
	&= \left| \int_{\mathcal{O}} \left[ F_{n, t}(x(t)) - F_{t}(x(t)) \right] w_{\phi,n}(t \mid \mathcal{O}) \, d t + \int_{\mathcal{O}} F_{t}(x(t)) \left[ w_{\phi,n}(t \mid \mathcal{O}) - w_{\phi}(t \mid \mathcal{O}) \right] \, d t \right| \nonumber \\
	&\leq \left| \int_{\mathcal{O}} \left[ F_{n, t}(x(t)) - F_{t}(x(t)) \right] w_{\phi,n}(t \mid \mathcal{O}) \, d t \right| + \left| \int_{\mathcal{O}} F_{t}(x(t)) \left[ w_{\phi,n}(t \mid \mathcal{O}) - w_{\phi}(t \mid \mathcal{O}) \right] \, d t \right|.
\end{align}

For the first part of ~\ref{eq:Kn_diff}, we claim that $0 \leq w_{\phi,n}(t \mid \mathcal{O}) \leq C$, where $C$ is a constant, for all $n$ and $t$. Bigger than 0 is by definition, we only need to prove $w_{\phi,n}(t \mid \mathcal{O})$ is bounded and the bound is independent with n and t. 

Note that $q_{n}(t)\in(0,1]$ and by definition of $\phi$ (continuous, bounded and positive), we can assume $\phi((0,1])\in(c,C]\quad,\forall t\in[0,1]$, where c is bigger than 0.Thus $ w_{\phi,n}(t \mid \mathcal{O}) = \frac{\phi(q_n(t))}{\int_{\mathcal{O}} \phi(q_n(t)) \, d t}\leq \frac{C}{c m(\mathcal{O})} := C$, where C is independent with $n$ and $t$, and WLOG we assume $m(\mathcal{O})>0$. Let
\[
w(t) = \begin{cases}
	w_{\phi,n}(t \mid \mathcal{O}) & \text{if } t \in \mathcal{O} \\
	0 & \text{otherwise}
\end{cases}
\]

Then we use this $w(t)$ in Theorem 1, we can see clearly that $C=1$. As a result we have 
$$
\sup _{x \in \operatorname{Lip}[0,1]} \left| \int_{\mathcal{O}} \left[ F_{n, t}(x(t)) - F_{t}(x(t)) \right] w_{\phi,n}(t \mid \mathcal{O}) \, d t \right| \xrightarrow{\text{a.s.}}0 
$$

For the second part of ~\ref{eq:Kn_diff}, we first consider
\[
\left| \int_{\mathcal{O}} \left[ w_{\phi}(t \mid \mathcal{O}) - w_{\phi,n}(t \mid \mathcal{O}) \right] \, d t \right|.
\]

By definition, $q_{n}(t) = \frac{\# I(t)}{n}$, where $I(t) := \{1 \leq i \leq n : t \in \mathcal{O}_{i}\}$ and $Q(t) = \mathbb{P}(\mathcal{O} \ni t)$. Thus, by the Strong Law of Large Numbers, we have
\begin{align}
	q_{n}(t) = \frac{\# I(t)}{n} = \frac{\sum_{i=1}^{n} \mathbf{1}_{\{ t \in \mathcal{O}_{i} \}}}{n}  \xrightarrow{\text{a.s.}} \mathbb{E}[\mathbf{1}_{\{ t \in \mathcal{O}_{i} \}}] = \mathbb{P}(\mathcal{O} \ni t) = Q(t) \nonumber
\end{align}

That means, as $n \rightarrow \infty$, by Lemma 4, and considering $\mathcal{O} \subset [0,1]$, we have $\int_{\mathcal{O}} | q_{n}(t) - Q(t) |\leq \int_{0}^{1} | q_{n}(t) - Q(t) |$, finally we get
\begin{align}
	\int_{\mathcal{O}} | q_{n}(t) - Q(t) | \, d t \xrightarrow{\text{a.s.}} 0.\nonumber
\end{align}

Now we can see $\int_{\mathcal{O}} q_{n}(t) \, d t \xrightarrow{\text{a.s.}} \int_{\mathcal{O}} Q(t) \, d t$, since $|\int_{\mathcal{O}}  q_{n}(t) - Q(t)  \, d t| \leq \int_{\mathcal{O}} | q_{n}(t) - Q(t) | \, d t$.

Noting that the $\int_{\mathcal{O}} q_{n}(t), \int_{\mathcal{O}} Q(t)$ is always positive by definition. And now we have $q_{n}(t)\xrightarrow{\text{a.s.}} Q(t)$ and $\int_{\mathcal{O}} q_{n}(t) \, d t \xrightarrow{\text{a.s.}} \int_{\mathcal{O}} Q(t) \, d t$. By Lemma 1, we get $w_{n}(t \mid \mathcal{O}) = \frac{q_n(t)}{\int_{\mathcal{O}} q_n(t) \, d t} \xrightarrow{\text{a.s.}} w(t \mid \mathcal{O}) = \frac{Q(t)}{\int_{\mathcal{O}} Q(t) \, d t}$. Noting that $\phi$ is a bounded continuous function, for $\phi$ not being the identity function, by Lemma 5, we have $\phi(q_{n}(t))\xrightarrow{\text{a.s.}} \phi(Q(t))$. Then replace $q_{n}$ and $Q$ with $\phi(q_{n})$ and $\phi(Q)$ respectively, we can get $\int_{\mathcal{O}} \phi(q_{n}(t)) \, d t \xrightarrow{\text{a.s.}} \int_{\mathcal{O}} \phi(Q(t)) \, d t$

following the same argument, we get
\[
w_{\phi,n}(t \mid \mathcal{O}) = \frac{\phi(q_n(t))}{\int_{\mathcal{O}} \phi(q_n(t)) \, d t} \xrightarrow{\text{a.s.}} w_{\phi}(t \mid \mathcal{O}) = \frac{\phi(Q(t))}{\int_{\mathcal{O}} \phi(Q(t)) \, d t}.
\]

Furthermore, replacing $q_{n}$ and $Q$ with $w_{\phi,n}$ and $w_{\phi}$, by lemma 4, we can similarly obtain
\[
\int_{\mathcal{O}} | w_{\phi}(t \mid \mathcal{O}) - w_{\phi,n}(t \mid \mathcal{O}) | \, d t \xrightarrow{\text{a.s.}} 0.\]

Considering that $F_t(x(t))$ is bounded by $1$, the second part of ~\ref{eq:Kn_diff}, i.e.,
\begin{align}
	\left| \int_{\mathcal{O}} F_{t}(x(t)) \left( w_{\phi}(t \mid \mathcal{O}) - w_{\phi,n}(t \mid \mathcal{O}) \right) \, d t \right| &\leq  \int_{\mathcal{O}} \left|F_{t}(x(t)) \left( w_{\phi}(t \mid \mathcal{O}) - w_{\phi,n}(t \mid \mathcal{O}) \right)\right|  \, d t\nonumber\\
	&\leq  \int_{\mathcal{O}} \left|\left( w_{\phi}(t \mid \mathcal{O}) - w_{\phi,n}(t \mid \mathcal{O}) \right)\right|  \, d t\nonumber\\
	& \xrightarrow{\text{a.s.}} 0\nonumber
\end{align}

and this convergence is independent with $x$. Finally, we obtain
\[
\sup_{x \in \operatorname{Lip}[0,1]} \left| K_{n}(x) - K(x) \right| \xrightarrow{\text{a.s.}} 0 .
\]
Thus,
\[
\sup_{x \in \operatorname{Lip}[0,1]} \left| POIFD_{n}(x) - POIFD(x) \right| \xrightarrow{\text{a.s.}} 0.
\]

\subsection{Proof of Theorem 3}

\textbf{Proof:} We consider a slight modification of the proof of Theorem 3.2 in Fraiman and Muniz (2001) \upcite{FraimanMuniz2001}. We define $\hat{\rho}_{n, h}$ and $\rho_{n, h}$ as follows:
\[
\hat{\rho}_{n, h}(t) = \frac{1}{n} \sum_{i=1}^{n} \mathbf{1}_{[\beta, +\infty)} \left( POIFD_{n} \left( X_{i} \right) \right) \mathbf{1}_{X_{i}(t) \text{ observable}} h \left( X_{i}(t) \right),
\]
\[
\rho_{n, h}(t) = \frac{1}{n} \sum_{i=1}^{n} \mathbf{1}_{[\beta, +\infty)} \left( POIFD \left( X_{i} \right) \right) \mathbf{1}_{X_{i}(t) \text{ observable}} h \left( X_{i}(t) \right).
\]

When \( h(t) = t \), \( \widehat{\rho}_{n,h} \) is the numerator of \( \widehat{\mu}_n \), and when \( h(t) = 1 \), then \( \widehat{\rho}_{n,h} \) is the denominator of \( \widehat{\mu}_n \). The same holds for \( \rho_{n,h} \) and \( \mu_n \).  

Since by the Strong Law of Large Numbers,

$$\rho_{n,h} \xrightarrow{\text{a.s.}} E\left[\mathbf{1}_{[\beta, +\infty)} \left( POIFD \left( X \right) \right) \mathbf{1}_{X(t) \text{ observable}} h \left( X(t) \right)\right]$$

Thus by Lemma 1, we have  $\hat\mu\xrightarrow{\text{a.s.}}\mu$.
Then to show $\hat\mu_{n}\xrightarrow{\text{a.s.}}\mu$, it suffices to show $\hat\mu_{n}\xrightarrow{\text{a.s.}}\hat\mu$. Also by Lemma 1,  it suffices to show that
\[
\left| \hat{\rho}_{n, h} - \rho_{n, h} \right| \xrightarrow{\text{a.s.}} 0.
\]

We define 
$$
S_{n} := \sup_{x \in \operatorname{Lip}[0,1]} \left| POIFD_{n}(x) - POIFD(x) \right|
$$

From Theorem 2, we know that $S_{n}\xrightarrow{\text{a.s.}} 0.$ For any $\delta > 0$,

\begin{align}
	&\left|\hat{\rho}_{n, h}-\rho_{n, h}\right|\nonumber \\
	& =  \left|\frac{1}{n} \sum_{i=1}^{n} \mathbf{1}_{[\beta, +\infty)} (\left( POIFD_{n} \left( X_{i} \right) \right)-\left( POIFD \left( X_{i} \right) \right)) \mathbf{1}_{X_{i}(t) \text{ observable}} h \left( X_{i}(t) \right)\right|\nonumber\\
	& =  \frac{1}{n} \sum_{i=1}^{n} \left|\mathbf{1}_{[\beta, +\infty)} (\left( POIFD_{n} \left( X_{i} \right) \right)-\left( POIFD \left( X_{i} \right) \right))\right| \left|\mathbf{1}_{X_{i}(t) \text{ observable}}\right|\left| h \left( X_{i}(t) \right)\right|\nonumber\\
	&\leq \frac{1}{n} \sum_{i=1}^{n}\left|h\left(X_{i}\right)\right|\left|\mathbf{1}_{X_{i}(t)\text{ is observed}}\right|\left|\mathbf{1}_{[\beta,+\infty)}\left(POIFD\left(X_{i}\right)+\delta\right)-\mathbf{1}_{[\beta,+\infty)}\left(POIFD\left(X_{i}\right)-\delta\right)\right| \mathbf{1}_{\left\{S_{n} \leq \delta\right\}} \nonumber\\
	&+\frac{1}{n} \sum_{i=1}^{n}\left|h\left(X_{i}\right)\right|\left|\mathbf{1}_{X_{i}(t)\text{ is observed}}\right|\left|\mathbf{1}_{[\beta,+\infty)}\left(POIFD\left(X_{i}\right)+S_{n}\right)-\mathbf{1}_{[\beta,+\infty)}\left(POIFD\left(X_{i}\right)-S_{n}\right)\right| \mathbf{1}_{\left\{S_{n} \geq \delta\right\}} \nonumber\\
	&\leq \frac{1}{n} \sum_{i=1}^{n}\left|h\left(X_{i}\right)\right|\left|\mathbf{1}_{[\beta,+\infty)}\left(POIFD\left(X_{i}\right)+\delta\right)-\mathbf{1}_{[\beta,+\infty)}\left(POIFD\left(X_{i}\right)-\delta\right)\right| \mathbf{1}_{\left\{S_{n} \leq \delta\right\}} \nonumber\\
	&+\frac{1}{n} \sum_{i=1}^{n}\left|h\left(X_{i}\right)\right|\left|\mathbf{1}_{[\beta,+\infty)}\left(POIFD\left(X_{i}\right)+S_{n}\right)-\mathbf{1}_{[\beta,+\infty)}\left(POIFD\left(X_{i}\right)-S_{n}\right)\right| \mathbf{1}_{\left\{S_{n} \geq \delta\right\}} \nonumber\\
	&\leq \frac{1}{n} \sum_{i=1}^{n}\left|h\left(X_{i}\right)\right|\left|\mathbf{1}_{[\beta,+\infty)}\left(POIFD\left(X_{i}\right)+\delta\right)-\mathbf{1}_{[\beta,+\infty)}\left(POIFD\left(X_{i}\right)-\delta\right)\right| \mathbf{1}_{\left\{S_{n} \leq \delta\right\}} \nonumber\\
	&+\frac{1}{n} \sum_{i=1}^{n}\left|h\left(X_{i}\right)\right| \mathbf{1}_{\left\{S_{n} \geq \delta\right\}} \nonumber\\
	&\xrightarrow{\text{a.s.}}E\left[\left|h\left(X\right)\right|\left|\mathbf{1}_{[\beta,+\infty)}\left(POIFD\left(X\right)+\delta\right)-\mathbf{1}_{\left[\beta,+\infty\right)}\left(POIFD\left(X\right)-\delta\right)\right|\right]\quad\text{By SLLN, as n $\rightarrow\infty$}.\nonumber\\
	&\text{ Since $E\left[h(X)\right]<\infty$, and indicator function is bounded, by Dominated Convergence Theorem,}\nonumber\\
	&\xrightarrow{\text{a.s.}} 0 \quad\text{as $\delta\rightarrow0$}. \nonumber
\end{align} 

Up to this point, we have completed all the proofs in the theory section.
\newpage
\section{Figures and Tables}
\begin{figure}
	\centering
	\subfigure[Partially observable functions]{
		\begin{minipage}[t]{0.4\textwidth}
			\centering
			\includegraphics[scale=0.45]{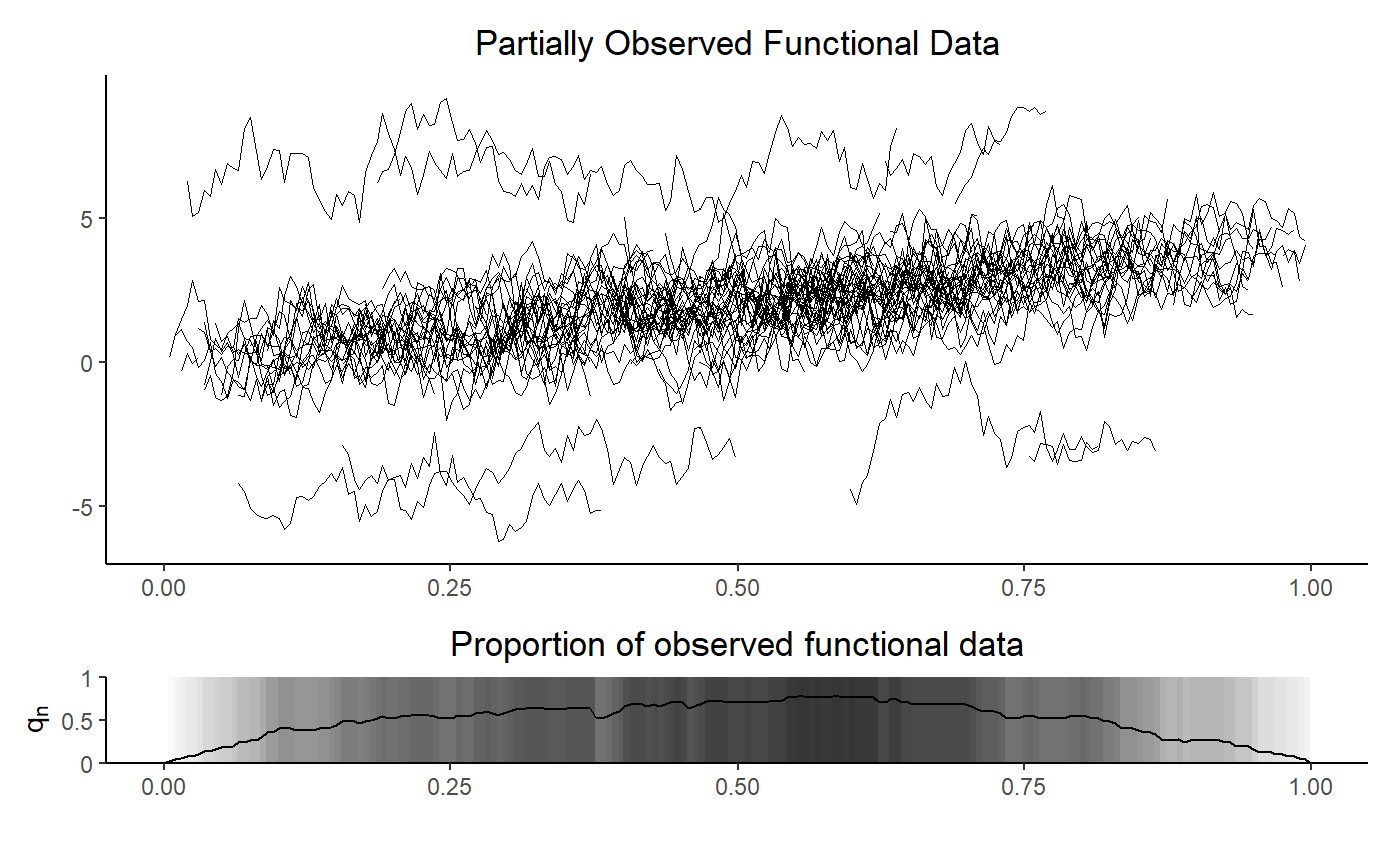}
		\end{minipage}
	}
	\hfill
	\subfigure[Trimming results]{
		\begin{minipage}[t]{0.4\textwidth}
			\centering
			\includegraphics[scale=0.45]{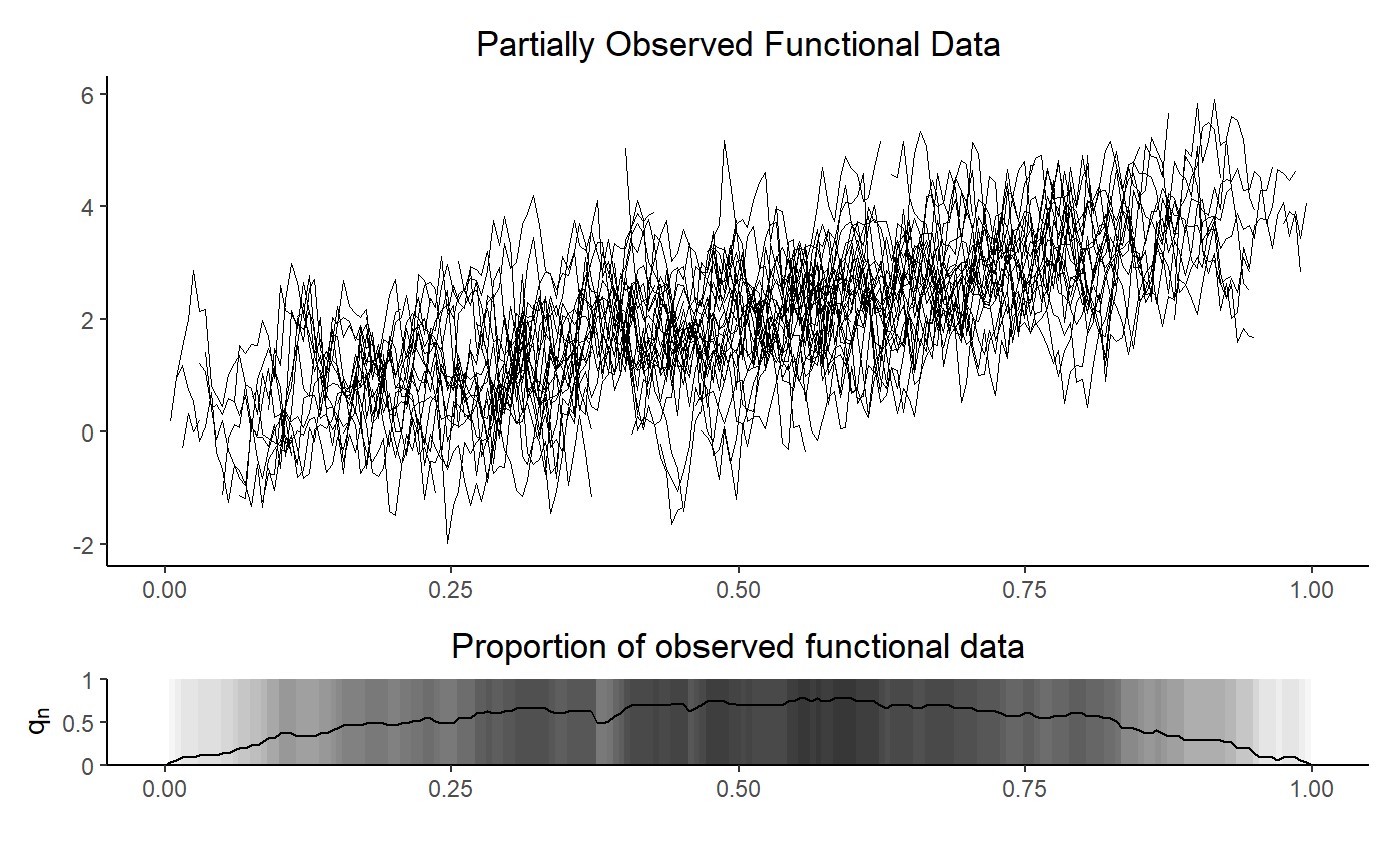}
		\end{minipage}
	}
	\centering
	\caption{Symmetric complete contamination model}
	\label{figure1}
\end{figure}

\begin{figure}
	\centering
	\subfigure[Partially observable functions]{
		\begin{minipage}[t]{0.4\textwidth}
			\centering
			\includegraphics[scale=0.45]{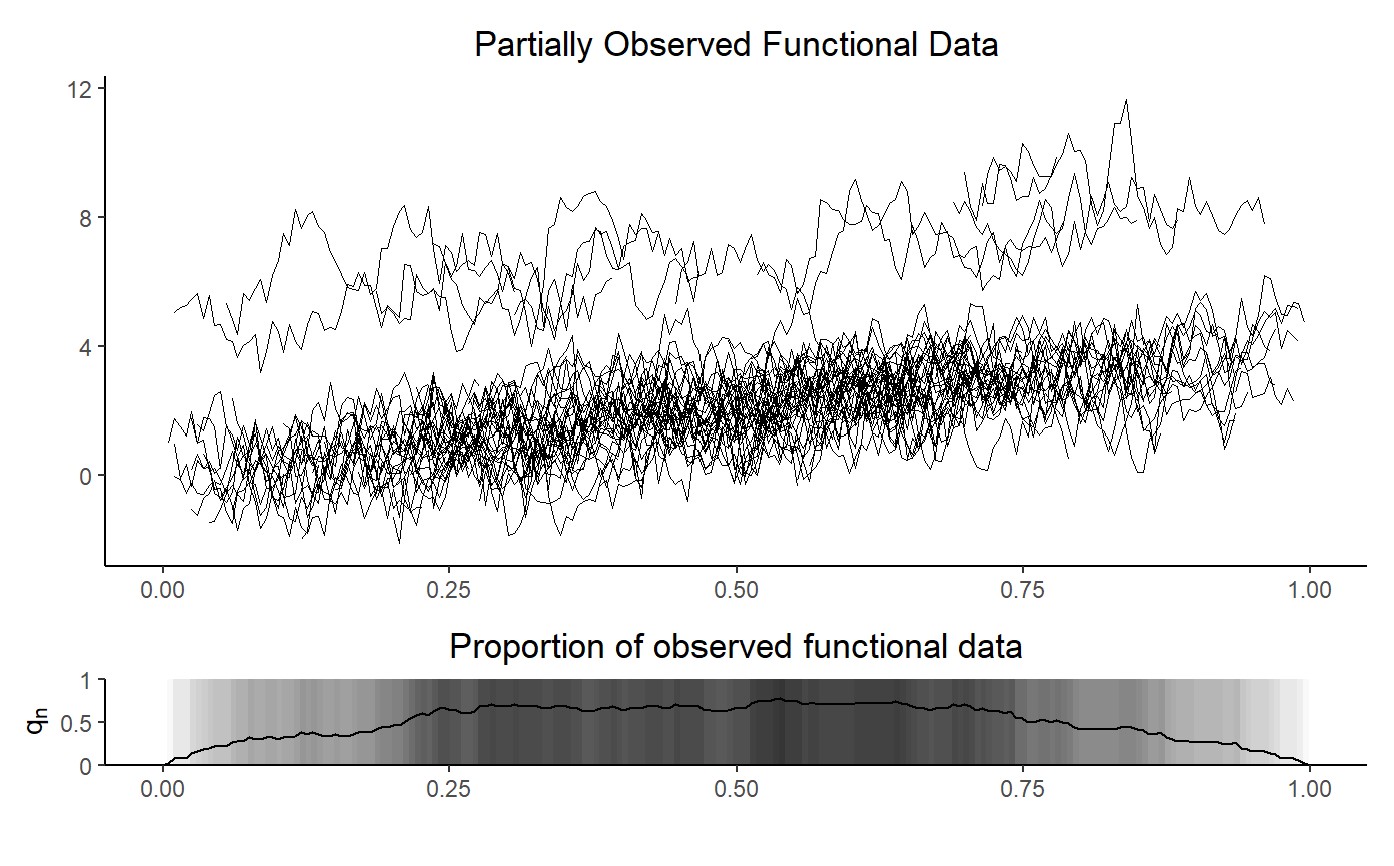}
		\end{minipage}
	}
	\hfill
	\subfigure[Trimming results]{
		\begin{minipage}[t]{0.4\textwidth}
			\centering
			\includegraphics[scale=0.45]{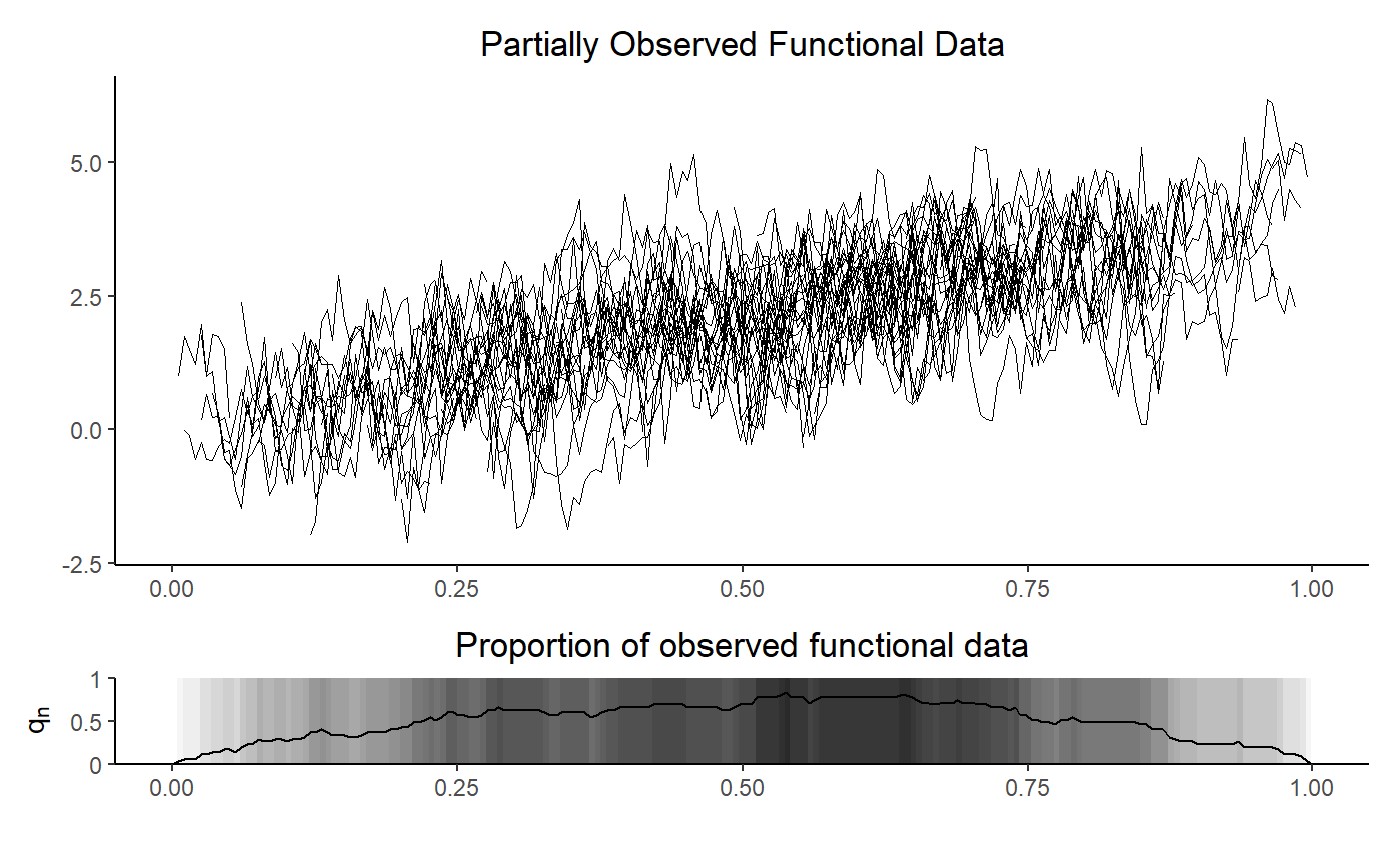}
		\end{minipage}
	}
	\centering
	\caption{Asymmetric complete contamination model}
	\label{figure2}
\end{figure}

\begin{figure}
	\centering
	\subfigure[Partially observable functions]{
		\begin{minipage}[t]{0.4\textwidth}
			\centering
			\includegraphics[scale=0.45]{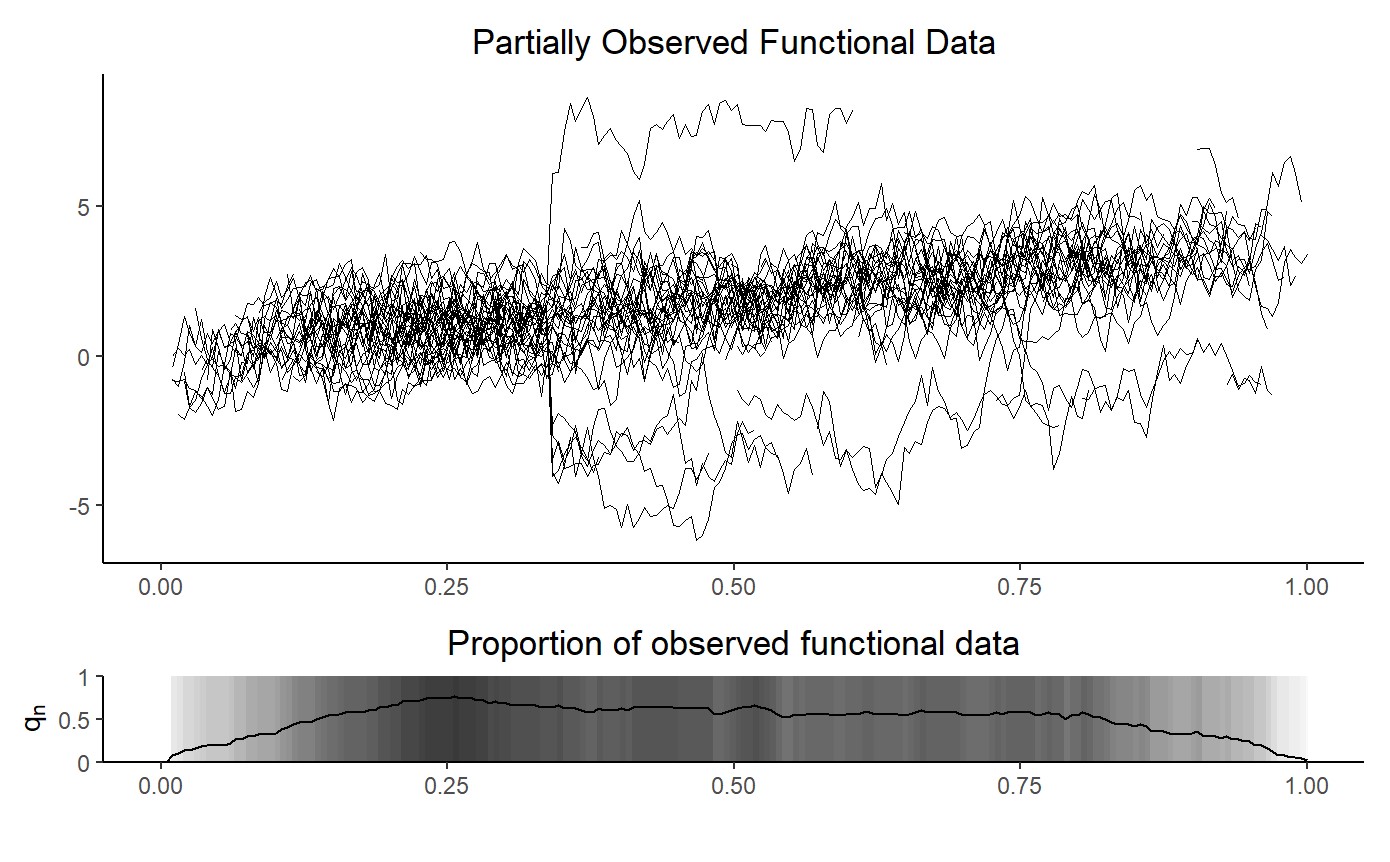}
		\end{minipage}
	}
	\hfill
	\subfigure[Trimming results]{
		\begin{minipage}[t]{0.4\textwidth}
			\centering
			\includegraphics[scale=0.45]{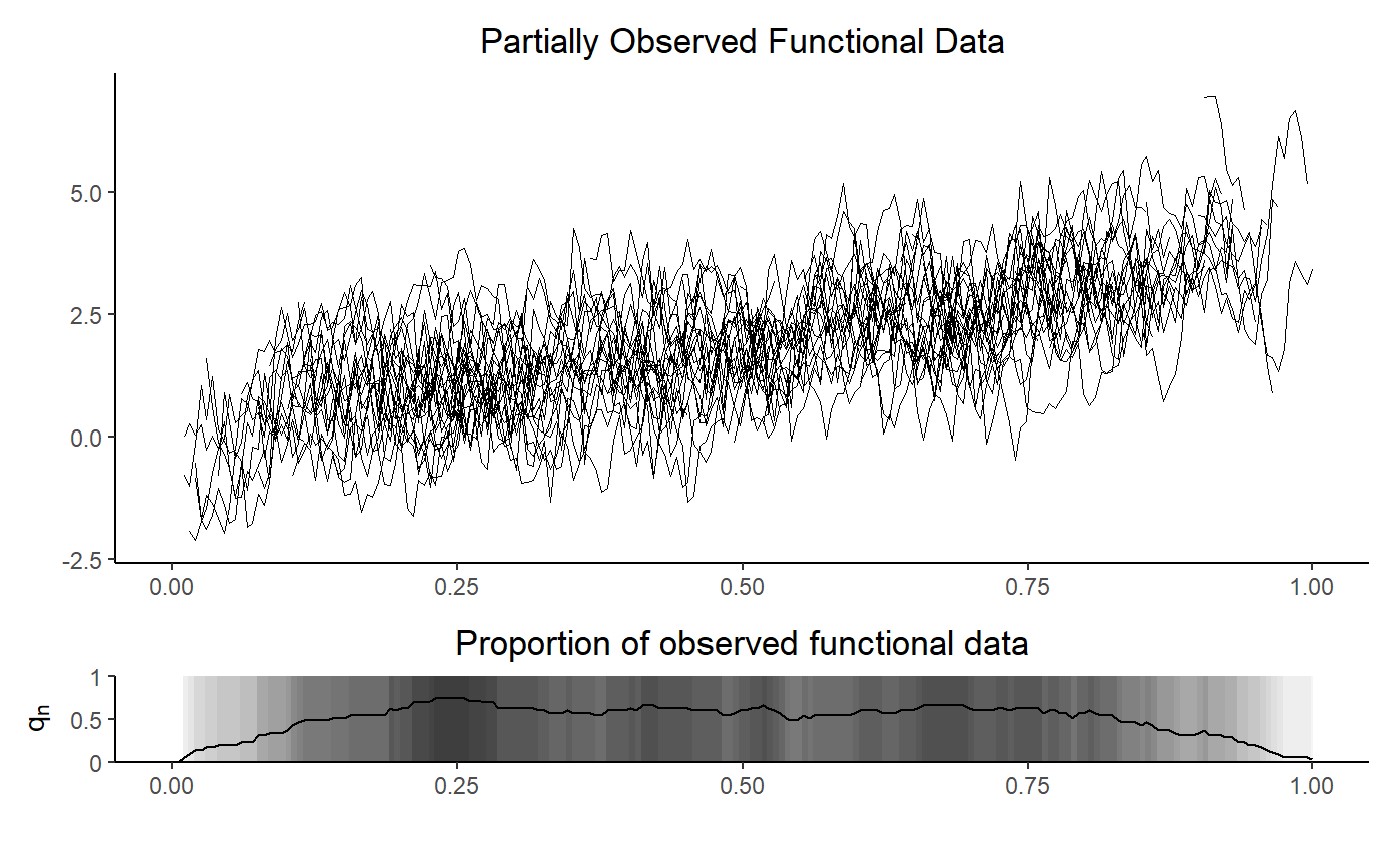}
		\end{minipage}
	}
	\centering
	\caption{Partial contamination model}
	\label{figure3}
\end{figure}
\begin{table}[ht]
	\centering
	\caption{$\alpha = 0.2$, observation proportion is 0.5}
	\resizebox{\textwidth}{!}{
		\begin{tabular}{rrrrrlrrrrrrr}
			\hline
			len & p & q & M & alpha & pollution\_type & observability & E & E\_trim & sd\_squared\_differences & sd\_squared\_differences\_trim & M & M\_trim \\ 
			\hline
			200.00 & 50.00 & 0.10 & 25.00 & 0.20 & symmetric & 0.50 & 7.14 & 0.51 & 7.93 & 1.34 & 3.70 & 0.08 \\ 
			200.00 & 50.00 & 0.10 & 25.00 & 0.20 & asymmetric & 0.50 & 9.59 & 0.08 & 4.91 & 0.02 & 8.85 & 0.07 \\ 
			200.00 & 50.00 & 0.10 & 25.00 & 0.20 & partial & 0.50 & 3.01 & 0.66 & 6.34 & 1.73 & 0.98 & 0.09 \\ 
			200.00 & 50.00 & 0.10 & 5.00 & 0.20 & symmetric & 0.50 & 0.24 & 0.09 & 0.15 & 0.04 & 0.16 & 0.08 \\ 
			200.00 & 50.00 & 0.10 & 5.00 & 0.20 & asymmetric & 0.50 & 0.43 & 0.10 & 0.26 & 0.04 & 0.42 & 0.08 \\ 
			200.00 & 50.00 & 0.10 & 5.00 & 0.20 & partial & 0.50 & 0.16 & 0.09 & 0.10 & 0.04 & 0.12 & 0.07 \\ 
			200.00 & 80.00 & 0.10 & 25.00 & 0.20 & symmetric & 0.50 & 2.03 & 0.05 & 1.57 & 0.01 & 1.47 & 0.06 \\ 
			200.00 & 80.00 & 0.10 & 25.00 & 0.20 & asymmetric & 0.50 & 8.26 & 0.06 & 6.38 & 0.02 & 6.63 & 0.06 \\ 
			200.00 & 80.00 & 0.10 & 25.00 & 0.20 & partial & 0.50 & 1.01 & 0.29 & 0.75 & 0.40 & 0.89 & 0.05 \\ 
			200.00 & 80.00 & 0.10 & 5.00 & 0.20 & symmetric & 0.50 & 0.12 & 0.05 & 0.07 & 0.01 & 0.10 & 0.05 \\ 
			200.00 & 80.00 & 0.10 & 5.00 & 0.20 & asymmetric & 0.50 & 0.42 & 0.06 & 0.27 & 0.02 & 0.36 & 0.06 \\ 
			200.00 & 80.00 & 0.10 & 5.00 & 0.20 & partial & 0.50 & 0.09 & 0.08 & 0.04 & 0.05 & 0.08 & 0.07 \\ 
			\hline
	\end{tabular}}
	
	\label{tab:tab1}
\end{table}

\begin{table}[ht]
	\centering
	\caption{$\alpha = 0.3$, observation proportion is 0.5}
	\resizebox{\textwidth}{!}{
		\begin{tabular}{rrrrrlrrrrrrr}
			\hline
			len & p & q & M & alpha & pollution\_type & observability & E & E\_trim & sd\_squared\_differences & sd\_squared\_differences\_trim & M & M\_trim \\ 
			\hline
			200.00 & 50.00 & 0.10 & 25.00 & 0.30 & symmetric & 0.50 & 2.81 & 0.08 & 2.00 & 0.01 & 1.83 & 0.08 \\ 
			200.00 & 50.00 & 0.10 & 25.00 & 0.30 & asymmetric & 0.50 & 8.10 & 0.10 & 8.23 & 0.03 & 4.57 & 0.09 \\ 
			200.00 & 50.00 & 0.10 & 25.00 & 0.30 & partial & 0.50 & 2.33 & 1.04 & 1.21 & 1.94 & 2.08 & 0.26 \\ 
			200.00 & 50.00 & 0.10 & 5.00 & 0.30 & symmetric & 0.50 & 0.33 & 0.09 & 0.30 & 0.01 & 0.23 & 0.09 \\ 
			200.00 & 50.00 & 0.10 & 5.00 & 0.30 & asymmetric & 0.50 & 0.43 & 0.10 & 0.17 & 0.04 & 0.44 & 0.09 \\ 
			200.00 & 50.00 & 0.10 & 5.00 & 0.30 & partial & 0.50 & 0.11 & 0.09 & 0.06 & 0.02 & 0.11 & 0.09 \\ 
			200.00 & 80.00 & 0.10 & 25.00 & 0.30 & symmetric & 0.50 & 2.12 & 0.07 & 1.55 & 0.03 & 1.83 & 0.06 \\ 
			200.00 & 80.00 & 0.10 & 25.00 & 0.30 & asymmetric & 0.50 & 8.26 & 0.07 & 3.55 & 0.02 & 7.86 & 0.06 \\ 
			200.00 & 80.00 & 0.10 & 25.00 & 0.30 & partial & 0.50 & 1.70 & 0.18 & 1.73 & 0.18 & 0.92 & 0.08 \\ 
			200.00 & 80.00 & 0.10 & 5.00 & 0.30 & symmetric & 0.50 & 0.15 & 0.08 & 0.07 & 0.04 & 0.15 & 0.07 \\ 
			200.00 & 80.00 & 0.10 & 5.00 & 0.30 & asymmetric & 0.50 & 0.42 & 0.06 & 0.26 & 0.02 & 0.30 & 0.06 \\ 
			200.00 & 80.00 & 0.10 & 5.00 & 0.30 & partial & 0.50 & 0.07 & 0.06 & 0.03 & 0.01 & 0.07 & 0.05 \\ 
			\hline
	\end{tabular}}
	
	\label{tab:tab2}
\end{table}

\begin{table}[ht]
	\centering
	\caption{$\alpha = 0.2$, observation proportion is 0.9}
	\resizebox{\textwidth}{!}{
		\begin{tabular}{rrrrrlrrrrrrr}
			\hline
			len & p & q & M & alpha & pollution\_type & observability & E & E\_trim & sd\_squared\_differences & sd\_squared\_differences\_trim & M & M\_trim \\ 
			\hline
			200.00 & 50.00 & 0.10 & 25.00 & 0.20 & symmetric & 0.90 & 2.36 & 0.03 & 2.25 & 0.01 & 1.51 & 0.03 \\ 
			200.00 & 50.00 & 0.10 & 25.00 & 0.20 & asymmetric & 0.90 & 14.07 & 0.08 & 9.58 & 0.13 & 13.53 & 0.04 \\ 
			200.00 & 50.00 & 0.10 & 25.00 & 0.20 & partial & 0.90 & 1.07 & 0.11 & 0.68 & 0.17 & 0.98 & 0.03 \\ 
			200.00 & 50.00 & 0.10 & 5.00 & 0.20 & symmetric & 0.90 & 0.07 & 0.03 & 0.03 & 0.01 & 0.07 & 0.03 \\ 
			200.00 & 50.00 & 0.10 & 5.00 & 0.20 & asymmetric & 0.90 & 0.35 & 0.03 & 0.33 & 0.01 & 0.19 & 0.03 \\ 
			200.00 & 50.00 & 0.10 & 5.00 & 0.20 & partial & 0.90 & 0.03 & 0.04 & 0.01 & 0.01 & 0.03 & 0.04 \\ 
			200.00 & 80.00 & 0.10 & 25.00 & 0.20 & symmetric & 0.90 & 1.53 & 0.02 & 1.66 & 0.00 & 0.74 & 0.02 \\ 
			200.00 & 80.00 & 0.10 & 25.00 & 0.20 & asymmetric & 0.90 & 5.33 & 0.02 & 2.95 & 0.00 & 4.94 & 0.02 \\ 
			200.00 & 80.00 & 0.10 & 25.00 & 0.20 & partial & 0.90 & 0.46 & 0.04 & 0.63 & 0.05 & 0.18 & 0.02 \\ 
			200.00 & 80.00 & 0.10 & 5.00 & 0.20 & symmetric & 0.90 & 0.06 & 0.02 & 0.04 & 0.00 & 0.04 & 0.02 \\ 
			200.00 & 80.00 & 0.10 & 5.00 & 0.20 & asymmetric & 0.90 & 0.29 & 0.02 & 0.21 & 0.00 & 0.20 & 0.02 \\ 
			200.00 & 80.00 & 0.10 & 5.00 & 0.20 & partial & 0.90 & 0.03 & 0.02 & 0.03 & 0.00 & 0.02 & 0.02 \\ 
			\hline
	\end{tabular}}
	
	\label{tab:tab3}
\end{table}

\begin{table}[ht]
	\centering
	\caption{$\alpha = 0.3$, observation proportion is 0.9}
	\resizebox{\textwidth}{!}{
		\begin{tabular}{rrrrrlrrrrrrr}
			\hline
			len & p & q & M & alpha & pollution\_type & observability & E & E\_trim & sd\_squared\_differences & sd\_squared\_differences\_trim & M & M\_trim \\ 
			\hline
			200.00 & 50.00 & 0.10 & 25.00 & 0.30 & symmetric & 0.90 & 1.24 & 0.04 & 0.98 & 0.01 & 0.87 & 0.04 \\ 
			200.00 & 50.00 & 0.10 & 25.00 & 0.30 & asymmetric & 0.90 & 7.75 & 0.04 & 5.43 & 0.01 & 7.73 & 0.04 \\ 
			200.00 & 50.00 & 0.10 & 25.00 & 0.30 & partial & 0.90 & 0.58 & 0.04 & 0.73 & 0.02 & 0.40 & 0.04 \\ 
			200.00 & 50.00 & 0.10 & 5.00 & 0.30 & symmetric & 0.90 & 0.11 & 0.05 & 0.10 & 0.01 & 0.07 & 0.05 \\ 
			200.00 & 50.00 & 0.10 & 5.00 & 0.30 & asymmetric & 0.90 & 0.53 & 0.04 & 0.21 & 0.01 & 0.45 & 0.04 \\ 
			200.00 & 50.00 & 0.10 & 5.00 & 0.30 & partial & 0.90 & 0.06 & 0.04 & 0.05 & 0.00 & 0.04 & 0.04 \\ 
			200.00 & 80.00 & 0.10 & 25.00 & 0.30 & symmetric & 0.90 & 1.44 & 0.02 & 1.23 & 0.00 & 1.03 & 0.03 \\ 
			200.00 & 80.00 & 0.10 & 25.00 & 0.30 & asymmetric & 0.90 & 8.09 & 0.02 & 3.76 & 0.00 & 7.34 & 0.02 \\ 
			200.00 & 80.00 & 0.10 & 25.00 & 0.30 & partial & 0.90 & 0.37 & 0.04 & 0.54 & 0.06 & 0.16 & 0.02 \\ 
			200.00 & 80.00 & 0.10 & 5.00 & 0.30 & symmetric & 0.90 & 0.06 & 0.02 & 0.04 & 0.01 & 0.06 & 0.02 \\ 
			200.00 & 80.00 & 0.10 & 5.00 & 0.30 & asymmetric & 0.90 & 0.32 & 0.03 & 0.27 & 0.00 & 0.23 & 0.03 \\ 
			200.00 & 80.00 & 0.10 & 5.00 & 0.30 & partial & 0.90 & 0.03 & 0.02 & 0.03 & 0.01 & 0.03 & 0.02 \\ 
			\hline
	\end{tabular}}
	
	\label{tab:tab4}
\end{table}

\end{document}